\documentclass[11pt]{amsart}
\usepackage[margin=1in]{geometry}
\usepackage{amssymb, mathtools} 
\usepackage[T1]{fontenc}
\usepackage{lmodern, microtype}
\usepackage[hidelinks]{hyperref}
\newtheorem{theorem}{Theorem}[section]
\newtheorem{lemma}[theorem]{Lemma}
\newtheorem{proposition}[theorem]{Proposition}
\newtheorem{corollary}[theorem]{Corollary}
\theoremstyle{definition}
\newtheorem{definition}[theorem]{Definition}
\theoremstyle{remark}
\newtheorem{remark}[theorem]{Remark}
\numberwithin{equation}{section}
\newcommand{\C}{\mathbb{C}}
\newcommand{\R}{\mathbb{R}}
\newcommand{\T}{\mathbb{T}}
\newcommand{\Prob}{\mathcal{P}}
\newcommand{\oB}{\overline{B}}
\DeclareMathOperator{\capc}{cap}
\DeclareMathOperator{\diam}{diam}

\DeclareMathOperator{\supp}{supp}

\title[Hausdorff continuity of logarithmic and Riesz capacities]{Hausdorff continuity of logarithmic \\and Riesz capacities}

\author{Sergei Kalmykov}
\address{School of Mathematical Sciences, Shanghai Jiao Tong University, 800 Dongchuan RD, Shanghai 200240, China \\and
\\
Department of Mathematics and Computer Sciences
Saint Petersburg State University, Saint Petersburg, 14 line of Vasilievsky island, 29, 199178, Russia
\\ and
\\
Khabarovsk Division of the Institute for Applied Mathematics, Far Eastern Branch, Russian Academy of Sciences, Khabarovsk, Seryshev st., 60, 680038, Russia.
} 
\email{kalmykovsergei@sjtu.edu.cn}

\author{Leonid V. Kovalev}
\address{Mathematics Department, Syracuse University, Syracuse, NY 13244, USA}
\email{lvkovale@syr.edu}



\begin{document}

\begin{abstract}
Let $E_n$ be $\alpha_n$-uniformly perfect compact sets that converge to a compact set $E$ in the Hausdorff distance $d_n=d_H(E_n,E)$. We prove that the logarithmic capacity of $E_n$ converges to that of $E$ whenever $\log(1/\alpha_n)=O(\log(1/d_n))$ and that no factor  tending to infinity can be inserted in this bound. A common energy-transfer theorem also yields Hausdorff continuity of Riesz capacity under local capacity bounds, including bounds obtained from lower Hausdorff content density. Finally, we correct our earlier claims about the uniform convergence of Green functions. 
\end{abstract}

\maketitle

\section{Introduction}\label{sec-intro}

Logarithmic capacity is upper semicontinuous under Hausdorff convergence of
compact sets, but it is not continuous: finite sets can approximate a
continuum. Continuity under a thickness assumption was proved by
Aseev~\cite{Aseev} for the conformal capacity of condensers whose plates are
uniformly perfect with a common constant, and by Aseev and
Lazareva~\cite{AseevLazareva} for the reduced modulus and the transfinite
diameter. Our earlier paper~\cite{KK} allowed the constants of uniform perfectness to
tend to zero. The purpose of this paper is threefold: close the gap between necessary and sufficient conditions in~\cite{KK}, establish similar convergence results for Riesz capacities; and to correct an erroneous claim of the uniform convergence of Green function in~\cite[Theorem 3.4]{KK}. 

\begin{definition}
Let $0<\alpha<1$. A compact set $F\subset\R^m$ with more than one point is
\emph{$\alpha$-uniformly perfect} if
\begin{equation}\label{eq-up}
 F\cap\{x\colon \alpha r\le|x-a|\le r\}\ne\varnothing
 \qquad\text{for all } a\in F \text{ and } 0<r\le\diam F.
\end{equation}
\end{definition}

Let $E_n$ and $E$ be nonempty compact subsets of $\R^m$, $m\ge1$, such
that $E_n$ is $\alpha_n$-uniformly perfect and $d_n=d_H(E_n,E)\to0$, where $d_H$ is the Hausdorff distance. We write
\[
 A_n=\log\frac1{\alpha_n},\qquad t_n=\log\frac1{d_n},
\]
omitting the terms with $d_n=0$. For a compact set $F\subset\R^m$ let
$\Prob(F)$ be the set of Borel probability measures on $F$, and let
\begin{equation}\label{eq-logcap}
 I(\mu)=\iint\log\frac1{|x-y|}\,d\mu(x)\,d\mu(y),\qquad
 V(F)=\inf_{\mu\in\Prob(F)}I(\mu),\qquad
 \capc F=e^{-V(F)}.
\end{equation}
In the plane, $\capc F$ is the logarithmic capacity~\cite{Ransford}; in other
dimensions it is the logarithmic energy capacity considered, for example,
in~\cite{CL}. By \cite[Theorem 3.4 and Remark 3.6]{KK}, planar sets satisfy
$\capc E_n\to\capc E$ provided that
\begin{equation}\label{eq-KK}
 A_n\le\frac1{24}\,\frac{t_n}{\log t_n}
\end{equation}
for all sufficiently large $n$, while the example of equally spaced arcs in
\cite[Proposition 4.1]{KK} shows that boundedness of $d_nA_n$ is not
sufficient. The following two theorems close most of the gap between these
conditions.

\begin{theorem}\label{thm-cap}
If $A_n=O(t_n)$, then $\capc E_n\to\capc E$.
\end{theorem}

\begin{theorem}\label{thm-sharp}
Let $H\colon[1,\infty)\to(0,\infty)$ be a nondecreasing function such that
$H(t)\to\infty$ as $t\to\infty$. There exist a compact set $E\subset[0,1]$ with
$\capc E>0$ and compact sets $E_n\subset[0,1]$ as above such that
$A_n\le H(t_n)\,t_n$ for all $n$, but $\limsup\capc E_n<\capc E$.
\end{theorem}

Theorem~\ref{thm-cap} is proved in Section~\ref{sec-transfer} by
redistributing an energy-minimizing measure on $E$ over small pieces of $E_n$.
Uniform perfectness enters only through a lower bound for the capacity of these
pieces, which in the plane is due to Pommerenke~\cite{Pommerenke}. The examples
in Theorem~\ref{thm-sharp} are sparse Cantor sets, constructed in
Section~\ref{sec-sharp}. 

Corollary~\ref{cor-circle} shows that the condition $d_nA_n \to0$ (which is much weaker than $A_n=O(t_n)$ in Theorem~\ref{thm-cap}) is sufficient
for capacity convergence when the limit is a circle or a segment. Proposition~4.1  in~\cite{KK} shows that this condition cannot be relaxed to $d_nA_n=O(1)$. As an application, the Cantor-type set obtained by removing the middle $\epsilon_k$-part of every interval at the $k$th step is NED in the sense of Ahlfors and
Beurling~\cite{AB} whenever $\log(1/\epsilon_k)=o(2^k)$
(Corollary~\ref{cor-ned}); the condition in \cite[Theorem 4.2]{KK} was
$\log(1/\epsilon_k)\le Ck/\log k$ with $C<1/(24\log 2)$.

In the plane, suppose $\capc E>0$, let $\Omega$ be the unbounded component of $\widehat\C\setminus E$
and let $g_E$ be its Green function with pole at $\infty$, extended by zero to
$\C\setminus\Omega$. In~\cite[Theorem 3.4]{KK} we claimed that $g_{E_n}\to g_E$ uniformly and $\Omega$ is a regular domain for the Dirichlet problem. This claim is false; in fact, the regularity may fail regardless of how slowly $A_n$ grow.  

\begin{theorem}\label{thm-green}
Let $E_n,E\subset\C$ be as above, with $A_n=O(t_n)$ and $\capc E>0$. Then
$g_{E_n}\to g_E$ uniformly on $\C$ if and only if $\Omega$ is regular for the
Dirichlet problem. For every function $\varphi\colon[1,\infty)\to(0,\infty)$ with
$\varphi(t)\to\infty$ there exist such sets with $A_n\le\varphi(t_n)$ for all
$n$ for which $\Omega$ is not regular.
\end{theorem}

The source of the error in \cite{KK}, and what remains valid there, is
described in Remark~\ref{rem-KK}. In particular, the capacity conclusion of
\cite[Theorem 3.4]{KK} and its application \cite[Theorem 4.2]{KK} are correct,
and they are improved by Theorem~\ref{thm-cap} and Corollary~\ref{cor-ned}.
Totik~\cite{Totik} proved that harmonic measures converge weakly when both the
sets and their capacities converge; Lemma~\ref{lem-green} is a uniform
statement for Green functions when the limit is regular, to be compared with
\cite{KKGreen}, where kernel convergence and bounded connectivity are assumed
instead of capacity convergence.

The transfer argument uses the logarithmic kernel only through a few of its
general properties. For the Riesz kernels $|x-y|^{-s}$ with $0<s<m$, it gives
capacity convergence whenever the pieces $E_n\cap\oB(a,h)$ at some scale
$h\ge d_n$ have Riesz capacity comparable to that of a ball of radius $h$
(Section~\ref{sec-riesz}). Uniformly perfect sets are exactly the compact sets of
positive lower Assouad dimension \cite[Lemma 2.1]{KLV}; see also~\cite{RSV}.
For Riesz kernels, a sufficient thickness condition is lower Hausdorff
content density with an exponent greater than $s$
(Theorem~\ref{thm-lowerreg}). This is related to the capacity density
conditions studied in~\cite{AH,CILV}; our estimates allow the density
constants to depend on $n$ and require a bound only at the transfer scale. 

For background on energy minimization in arbitrary dimensions, we refer
to~\cite{BHS,Landkof,Mattila}. Clark and Laugesen~\cite{CL} study continuity
as the Riesz exponent varies on a fixed set; here the exponent is fixed and
the sets vary. The potential bound and the content-to-capacity estimate used
below are standard potential-theoretic tools. We record their short proofs
to keep the constants explicit and to treat all $0\le s<m$ uniformly.

\section{Energy transfer}\label{sec-transfer}

Let $B(a, r)$ and $\overline{B}(a, r)$ denote the open (respectively, closed) ball with center $a$ and radius $r$. 
To treat the logarithmic and Riesz kernels simultaneously, for any ambient dimension $m\ge1$ and any order of singularity $s\in [0, m)$, we let
\[
 k_0(x)=\log\frac1{|x|},\qquad k_s(x)=|x|^{-s}\quad(0<s<m).
\]
The value at the origin is $+\infty$. For finite positive Borel measures $\mu,\nu$ with compact support let
\[
 I_s(\mu,\nu)=\iint k_s(x-y)\,d\mu(x)\,d\nu(y),\qquad I_s(\mu)=I_s(\mu,\mu),
 \qquad U_s^\mu(x)=\int k_s(x-y)\,d\mu(y),
\]
and for a compact set $F$ let $V_s(F)=\inf_{\mu\in\Prob(F)}I_s(\mu)$. For $\mu\in\Prob(F)$ put
\[
 Q_\mu(r)=(\mu\times\mu)\{(x,y)\colon|x-y|\le r\}.
\]
A measure of finite energy has no atoms, and for $0<r<1$
\begin{equation}\label{eq-tail}
 k_s(r)\,Q_\mu(r)\le\iint_{|x-y|\le r}k_s(x-y)\,d\mu(x)\,d\mu(y)
 \to 0\quad(r\to0),
\end{equation}
where $k_s(r)$ is the value of $k_s$ on the sphere of radius $r$. Since the
kernel is lower semicontinuous and bounded below on bounded sets, weak
compactness of $\Prob(F)$ shows that $V_s(F)$ is attained when it is finite, by
an \emph{equilibrium measure}, and that
\begin{equation}\label{eq-usc}
 \liminf_{n\to\infty}V_s(F_n)\ge V_s(F)\qquad\text{whenever } d_H(F_n,F)\to0.
\end{equation}
This is the upper semicontinuity of capacity. The equilibrium measure is
unique when the energy is finite; see~\cite[Theorems 4.4.5 and 4.4.8]{BHS}
for the Riesz and logarithmic cases, respectively.

The next lemma combines the variational property of an equilibrium measure
with a weak maximum principle. For nonnegative radially decreasing kernels,
a general form of the latter is given in~\cite[Theorem 2.6.2]{AH}. The
elementary argument below also gives an additive bound for the logarithmic
kernel in every dimension.

\begin{lemma}\label{lem-maximum}
Let $F\subset\R^m$ be a compact set with $V_s(F)<\infty$, and let $\sigma$ be
an equilibrium measure of $F$. Then $U_s^\sigma\le V_s(F)$ on $\supp\sigma$, and
for all $x\in\R^m$
\[
 U_0^\sigma(x)\le V_0(F)+\log2,\qquad
 U_s^\sigma(x)\le 2^sV_s(F)\quad(0<s<m).
\]
\end{lemma}

\begin{proof}
Write $V=V_s(F)$ and $U=U_s^\sigma$. Given $\epsilon>0$, the set
$S=\{x\in F\colon U(x)\le V-\epsilon\}$ is closed because $U$ is lower
semicontinuous. If $\sigma(S)>0$, the measure $\tau=\sigma(S)^{-1}\sigma|_S$ has
finite energy and $I_s(\sigma,\tau)=\int U\,d\tau\le V-\epsilon$, hence
\[
 I_s\bigl((1-t)\sigma+t\tau\bigr)
 =V-2t\bigl(V-I_s(\sigma,\tau)\bigr)+t^2\bigl(V-2I_s(\sigma,\tau)+I_s(\tau)\bigr)<V
\]
for small $t>0$, contrary to the minimality of $\sigma$. Therefore
$U\ge V$ $\sigma$-almost everywhere, and since $\int U\,d\sigma=V$, equality
holds $\sigma$-almost everywhere. If $U(x_0)>V$ for some $x_0\in\supp\sigma$,
then $U>V$ on a neighborhood of $x_0$, which has positive $\sigma$-measure.
Thus $U\le V$ on $\supp\sigma$.

For $x\notin\supp\sigma$ let $z$ be a point of $\supp\sigma$ nearest to $x$.
For $y\in\supp\sigma$ we have $|z-y|\le|z-x|+|x-y|\le2|x-y|$, hence
$k_0(x-y)\le k_0(z-y)+\log 2$ and $k_s(x-y)\le2^sk_s(z-y)$ for $s>0$.
Integrating against $\sigma$ and using $U(z)\le V$ completes the proof.
\end{proof}

\begin{theorem}\label{thm-transfer}
Let $E_n,E\subset\R^m$ be nonempty compact sets and let $h_n\to0$ be positive
numbers with $d_H(E_n,E)\le h_n$. Suppose that
\begin{equation}\label{eq-W}
 W_n:=\sup_{a\in E_n}V_s\bigl(E_n\cap\oB(a,h_n)\bigr)<\infty .
\end{equation}
If $\mu\in\Prob(E)$ satisfies $I_s(\mu)<\infty$ and
$W_n\,Q_\mu(12h_n)\to0$, then there are $\mu_n\in\Prob(E_n)$ such that
$\mu_n\to\mu$ weakly and $\limsup_n I_s(\mu_n)\le I_s(\mu)$.
\end{theorem}

\begin{proof}
Let $\{b_j\}$ be a finite $h_n$-net in $E_n$. Since every point of $E$ lies
within $h_n$ of $E_n$, there is a Borel partition $\{P_j\}$ of $E$ such that
$|x-b_j|\le2h_n$ for $x\in P_j$. Let $\sigma_j$ be an equilibrium measure of
$F_j=E_n\cap\oB(b_j,h_n)$, which exists because $V_s(F_j)\le W_n$, and let
\[
 \mu_n=\sum_j\mu(P_j)\,\sigma_j\in\Prob(E_n).
\]
Then $I_s(\mu_n)=\iint K_n\,d\mu\,d\mu$, where
$K_n(x,y)=I_s(\sigma_i,\sigma_j)$ for $x\in P_i$, $y\in P_j$. By
Lemma~\ref{lem-maximum}, $K_n\le 2^sW_n+\log2$.

Every point of $\supp\sigma_i$ lies within $3h_n$ of every point of $P_i$.
This proves $\mu_n\to\mu$ weakly by uniform continuity of test functions
on a fixed compact set containing all the supports. It also gives
$K_n(x,y)\to k_s(x-y)$ for $x\ne y$. If $|x-y|>12h_n$, all distances
between $\supp\sigma_i$ and $\supp\sigma_j$ lie between $|x-y|/2$ and
$3|x-y|/2$, which gives $|K_n(x,y)|\le2^s|k_s(x-y)|+\log2$. The right-hand side
is $(\mu\times\mu)$-integrable, and the diagonal is $(\mu\times\mu)$-null. By
dominated convergence the integral of $K_n$ over $\{|x-y|>12h_n\}$ tends to
$I_s(\mu)$, while the integral over $\{|x-y|\le12h_n\}$ is at most
$(2^sW_n+\log 2)\,Q_\mu(12h_n)\to0$.
\end{proof}

\begin{corollary}\label{cor-transfer}
In the setting of Theorem~\ref{thm-transfer}, suppose $V_s(E)<\infty$ and
there are finite-energy measures $\mu_k\in\Prob(E)$ with
$I_s(\mu_k)\to V_s(E)$ such that
$W_nQ_{\mu_k}(12h_n)\to0$ as $n\to\infty$, for each $k$. Then
$V_s(E_n)\to V_s(E)$. This applies, in particular, when $V_s(E)<\infty$ and
the equilibrium measure of $E$ satisfies $W_nQ_\mu(12h_n)\to0$.
\end{corollary}

\begin{proof}
Theorem~\ref{thm-transfer} gives $\limsup_nV_s(E_n)\le I_s(\mu_k)$ for every
$k$, and~\eqref{eq-usc} gives the reverse inequality.
\end{proof}

The following estimate is due to Pommerenke~\cite[Theorem 1]{Pommerenke} in the plane. We include a construction that proves the same bound in every dimension.

\begin{lemma}\label{lem-localcap}
If $F\subset\R^m$ is $\alpha$-uniformly perfect, then
\[
 \capc\bigl(F\cap\oB(a,r)\bigr)\ge\frac{\alpha^2r}{32},
 \qquad a\in F,\quad 0<r\le\diam F.
\]
\end{lemma}

\begin{proof}
Let $\eta=\alpha/8$. If $b\in F$ and $0<\rho\le\diam F$, then~\eqref{eq-up}
with radius $\rho/2$ provides $b'\in F$ with $\alpha\rho/2\le|b-b'|\le\rho/2$.
The balls $\oB(b,\eta\rho)$ and $\oB(b',\eta\rho)$ are contained in
$\oB(b,\rho)$, and the distance between them is at least $\alpha\rho/4$.
Starting with $\oB(a,r)$ and repeating this construction, we obtain $2^k$
nested balls of radius $\eta^kr$ at level $k$, centered at points of $F$; their
limit set lies in $F\cap\oB(a,r)$. Let $\sigma$ be the probability measure that
assigns mass $2^{-k}$ to each ball of level $k$, and let $D_k$ be the set of
pairs $(x,y)$ that lie in a common ball of level $k$ but in different balls of
level $k+1$. Then
$(\sigma\times\sigma)(D_k)=2^k\bigl(2^{-2k}-2\cdot2^{-2k-2}\bigr)=2^{-k-1}$ and
$|x-y|\ge\alpha\eta^kr/4$ on $D_k$. The sets $D_k$ cover the complement of the
diagonal, which is $(\sigma\times\sigma)$-null. Therefore
\[
 I(\sigma)\le\sum_{k=0}^\infty2^{-k-1}
 \left(\log\frac4{\alpha r}+k\log\frac8\alpha\right)
 =\log\frac{32}{\alpha^2r}. \qedhere
\]
\end{proof}

\begin{proof}[Proof of Theorem~\ref{thm-cap}]
If $\capc E=0$, the claim follows from~\eqref{eq-usc}. Otherwise $E$ has more
than one point, so $d_n\le\diam E_n$ for large $n$, and
Lemma~\ref{lem-localcap} shows that~\eqref{eq-W} holds with $h_n=d_n$ and
\begin{equation}\label{eq-Wlog}
 W_n\le2A_n+t_n+\log32 .
\end{equation}
By~\eqref{eq-tail}, the equilibrium measure $\nu$ of $E$ satisfies
$t_nQ_\nu(12d_n)\to0$. Hence $W_nQ_\nu(12d_n)\to0$ when $A_n=O(t_n)$, and
Corollary~\ref{cor-transfer} completes the proof.
\end{proof}

\section{Circles, segments, and NED sets}\label{sec-circle}

When the limit set carries near-minimizing measures of controlled
concentration, Theorem~\ref{thm-transfer} gives conditions of a different type.

\begin{corollary}\label{cor-growth}
Suppose that $V(E)<\infty$ and that there are $\gamma>0$, measures
$\mu_k\in\Prob(E)$ with $I(\mu_k)\to V(E)$, and constants $M_k$ such that
$Q_{\mu_k}(r)\le M_kr^\gamma$ for $0<r\le1$. If $d_n^\gamma A_n\to0$, then
$\capc E_n\to\capc E$.
\end{corollary}

\begin{proof}
By~\eqref{eq-Wlog},
$W_nQ_{\mu_k}(12d_n)\le12^\gamma M_k(2A_n+t_n+\log32)\,d_n^\gamma\to0$, since
also $d_n^\gamma t_n\to0$. Apply Corollary~\ref{cor-transfer}.
\end{proof}

A growth bound $\mu(\oB(x,r))\le Mr^\gamma$ for $x\in E$ implies
$Q_\mu(r)\le Mr^\gamma$, which is how Corollary~\ref{cor-growth} is usually
applied.

\begin{corollary}\label{cor-circle}
If $E$ is a circle or a nondegenerate segment in $\R^m$, then
$d_nA_n\to0$ implies $\capc E_n\to\capc E$.
\end{corollary}

\begin{proof}
The energy~\eqref{eq-logcap} depends only on the distances between points, and
the hypothesis is invariant under scaling of all sets. So we may regard $E$ as
the unit circle or the segment $[0,1]$ in $\C$. The equilibrium measure of the
circle is normalized arclength, which satisfies $Q(r)\le r$. The equilibrium
measure of $[0,1]$ has density $(x(1-x))^{-1/2}/\pi$. Its normalized
restriction $\mu_\eta$ to $[\eta,1-\eta]$ has bounded density, so
$Q_{\mu_\eta}(r)\le C_\eta r$, and $I(\mu_\eta)\to V([0,1])$ as $\eta\to0$ by
dominated convergence. In both cases Corollary~\ref{cor-growth} applies with
$\gamma=1$.
\end{proof}


The condition in Corollary~\ref{cor-circle} cannot be replaced by boundedness
of $d_nA_n$, even with a small bound. For $c>0$ the sets
$E_n=\{z\in\T\colon|\arg z^n|\le e^{-cn}\}$ satisfy $d_H(E_n,\T)\le\pi/n$, they
are $\alpha_n$-uniformly perfect with $A_n\le cn+O(1)$, and
$\capc E_n=(\sin(e^{-cn}/2))^{1/n}\to e^{-c}$; see the proof of
\cite[Proposition 4.1]{KK}. Thus $d_nA_n\le\pi c+o(1)$, while
$\lim_n\capc E_n<\capc\T$. Theorem~\ref{thm-sharp} with $H(t)=t$ shows that for a general
limit set the condition $d_nA_n\to0$ is not sufficient. 

The following corollary refers to NED sets (negligible for extremal distances) introduced by Ahlfors and Beurling in~\cite{AB}. 

\begin{corollary}\label{cor-ned}
Let $0<\epsilon_k<1$, let $K_0=[0,1]$, and let $K_k$ be obtained from
$K_{k-1}$ by removing the open middle interval of relative length $\epsilon_k$
from each component. If $\log(1/\epsilon_k)=o(2^k)$, then
$K=\bigcap_kK_k$ satisfies $\capc([0,1]\setminus K)=\capc[0,1]$, where the left
side is the inner capacity. Consequently, $K$ is NED.
\end{corollary}

\begin{proof}
Let $E_n$ be the closure of $[0,1]\setminus K_n$. By the proof of
\cite[Theorem 4.2]{KK}, $E_n$ is $\alpha_n$-uniformly perfect with
$\alpha_n=\frac12\min_{k\le n}\epsilon_k$. The components of $K_n$ have length
$\ell_n=2^{-n}\prod_{k\le n}(1-\epsilon_k)$. Each point of such a component is
within $\ell_n/2$ of $E_n$, except in the two end components, where the
distance is at most $\ell_n$ and is attained at $0$ and $1$. Hence
$d_n=\ell_n\le2^{-n}$. Let $a_k=\log(1/\epsilon_k)$. Given $\delta>0$ we have
$a_k\le\delta2^k$ for $k\ge k_\delta$, hence
$\max_{k\le n}a_k\le\max_{k<k_\delta}a_k+\delta2^n$. Thus
$\max_{k\le n}a_k=o(2^n)$ and
\[
 d_nA_n\le2^{-n}\Bigl(\log2+\max_{k\le n}a_k\Bigr)\to0.
\]
By Corollary~\ref{cor-circle}, $\capc E_n\to\capc[0,1]$. Each $E_n$ is a finite
union of nondegenerate closed intervals whose interiors lie in
$[0,1]\setminus K$. Shrinking these intervals slightly produces compact subsets
of $[0,1]\setminus K$ that are $\alpha$-uniformly perfect with a fixed $\alpha$
and converge to $E_n$, so their capacities tend to $\capc E_n$ by
Theorem~\ref{thm-cap}. Therefore $\capc([0,1]\setminus K)=\capc[0,1]$. By a
theorem of Ahlfors and Beurling \cite[Theorem 14]{AB}, a compact subset $K$ of a segment $I$ is NED if and only if $\capc(I\setminus K)=\capc I$.
\end{proof}

\section{Sharpness of the capacity criterion}\label{sec-sharp}

The idea of the construction is to choose the limit set $E$ to be a
Cantor-type set which looks self-similar on most scales, except for a few
carefully selected ones. On the selected scales the subintervals of pre-Cantor
sets are made smaller than normal, but still not so small as to make $E$ a
polar set. The approximating sets $E_j$ are obtained by shrinking each
subinterval of such a pre-Cantor set toward its left endpoint: this reduces
their capacity well below $\capc E$, and comes at some cost to the uniform
perfectness constants.

\begin{proof}[Proof of Theorem~\ref{thm-sharp}]
Since $H(t)\to\infty$, we can choose positive integers $k_1<k_2<\cdots$ with
\begin{equation}\label{eq-kj}
 2^{k_j}\ge j^2\qquad\text{and}\qquad H\Bigl(\frac{2^{k_j}}{j^2}\Bigr)\ge7j^2.
\end{equation}
For integers $k\ge0$ let
\[
 T_k=k\log4+\sum_{j\colon k_j\le k}\frac{2^{k_j}}{j^2},\qquad \ell_k=e^{-T_k}.
\]
Then $\ell_0=1$ and $\ell_k\le\ell_{k-1}/4$. Let $C_0=[0,1]$, and obtain the
pre-Cantor set $C_k$ by keeping in each component of $C_{k-1}$ its leftmost
and rightmost subintervals of length $\ell_k$. The $2^k$ components of $C_k$
are the intervals of level $k$. Let $E=\bigcap_kC_k$, and let
$\mu\in\Prob(E)$ assign mass $2^{-k}$ to every interval of level $k$.

Let $D_k$, $k\ge1$, be the set of pairs $(x,y)\in E\times E$ that lie in a
common interval of level $k-1$ but in different intervals of level $k$. Then
$(\mu\times\mu)(D_k)=2^{k-1}\bigl(2^{-2k+2}-2\cdot2^{-2k}\bigr)=2^{-k}$, and
$|x-y|\ge\ell_{k-1}-2\ell_k\ge\ell_{k-1}/2$ on $D_k$. The sets $D_k$ cover the
complement of the diagonal, which is $(\mu\times\mu)$-null. Hence
\[
 I(\mu)\le\sum_{k\ge1}2^{-k}(T_{k-1}+\log2)
 =\log2+\frac12\sum_{k\ge0}2^{-k}T_k
 =\log8+\sum_{j\ge1}\frac1{j^2}<4,
\]
which shows that $\capc E>e^{-4}$.

Let $N_j=2^{k_j}$ and $r_j=e^{-6N_j}$. Since $2^{k_i}\le N_j$ for $i\le j$, we
have $N_j/j^2\le T_{k_j}\le k_j\log4+\frac{\pi^2}{6}N_j<6N_j$, so
$r_j<\ell_{k_j}$. Let $E_j$ be obtained from $C_{k_j}$ by replacing each
interval $[a,a+\ell_{k_j}]$ of level $k_j$ with $[a,a+r_j]$. Thus
$E_j\subset C_{k_j}\subset[0,1]$ is a union of $N_j$ disjoint intervals of
length $r_j$, and $\diam E_j\le1$. From a point of such an interval, all
distances in $[0,r_j/2]$ are attained within the same interval. Hence $E_j$ is
$\alpha_j$-uniformly perfect with $\alpha_j=r_j/2$.

Both $E$ and $E_j$ are contained in $C_{k_j}$ and meet every interval of level
$k_j$. Hence $d_j=d_H(E_j,E)\le\ell_{k_j}=e^{-T_{k_j}}$, so that
$t_j=\log(1/d_j)\ge T_{k_j}\ge N_j/j^2\ge1$, and by~\eqref{eq-kj}
\[
 A_j=6N_j+\log2\le7N_j\le H\Bigl(\frac{N_j}{j^2}\Bigr)\frac{N_j}{j^2}
 \le H(t_j)\,t_j.
\]

Finally, let $\nu\in\Prob(E_j)$ assign the masses $m_1,\dots,m_{N_j}$ to the
intervals of $E_j$. An interval of length $r_j$ has capacity $r_j/4$, so the
energy of $\nu$ within the $i$th interval is at least
$m_i^2\log(4/r_j)\ge6N_jm_i^2$. Since $\diam E_j\le1$, the kernel is
nonnegative on $E_j\times E_j$, and therefore
\[
 I(\nu)\ge6N_j\sum_im_i^2\ge6.
\]
Thus $\capc E_j\le e^{-6}<e^{-4}<\capc E$.
\end{proof}

\section{Green functions in the plane}\label{sec-green}

Let $E\subset\C$ be a compact set with $\capc E>0$, let $\Omega$ be the unbounded component of $\widehat\C\setminus E$, and let $\nu_E$ be the equilibrium measure of $E$. Since this section is concerned only with logarithmic capacity, we write $U^\nu =U_0^\nu$ for the energy of $\nu$ and $V(E)=V_0(E)$ for the infimal energy among $\nu\in \Prob(E)$. The Green function of $\Omega$ with pole at
$\infty$ is $g_E=V(E)-U^{\nu_E}$ on $\Omega$, and we set $g_E=0$ on
$\C\setminus\Omega$. When $\Omega$ is regular, $g_E$ is continuous on $\C$ and
the identity $g_E=V(E)-U^{\nu_E}$ holds on all of $\C$; in particular,
$g_E$ is subharmonic on $\C$~\cite[Section 4.4]{Ransford}. Every uniformly
perfect compact set has a regular outer domain; the quantitative
H\"older estimate in~\cite[Theorem 3.2]{KK} gives this directly.

\begin{lemma}\label{lem-green}
Let $E_n,E\subset\C$ be compact sets with regular outer domains. If
$d_H(E_n,E)\to0$ and $\capc E_n\to\capc E>0$, then $g_{E_n}\to g_E$
uniformly on $\C$.
\end{lemma}

\begin{proof}
Write $g_n=g_{E_n}$, $g=g_E$, and let $\nu_n,\nu$ be the equilibrium measures.
Every subsequential weak limit of $(\nu_n)$ is supported on $E$ and, by lower
semicontinuity of energy, has energy at most $\lim_n V(E_n)=V(E)$; by uniqueness
of the equilibrium measure it equals $\nu$. Thus $\nu_n\to\nu$ weakly.
For $z_n\to z$, weak convergence and lower semicontinuity of the kernel give
\[
 \liminf_{n\to\infty}U^{\nu_n}(z_n)\ge U^\nu(z).
\]
Indeed, this follows first for the continuous truncated kernels
$\min(k_0,M)$ and then by letting $M\to\infty$. Since $V(E_n)\to V(E)$,
we obtain $\limsup_n g_n(z_n)\le g(z)$. Continuity of $g$ and a compactness
argument therefore give
\[
 \limsup_{n\to\infty} \, \sup_K(g_n-g)\le0
 \qquad\text{for every compact }K\subset\C.
\]
On each component of $\C\setminus E_n$ the function $g-g_n$ is subharmonic, and on
the unbounded component it has a removable singularity at $\infty$. Its upper
limits at the boundary do not exceed $\sup_{E_n}g$. The maximum principle gives
$g-g_n\le\sup_{E_n}g$ on $\C$, and $\sup_{E_n}g\to0$ because $g$ vanishes
on $E$ and is uniformly continuous on a fixed compact neighborhood
containing the sets $E_n$. Thus $g_n\to g$ uniformly on compact sets.
Outside a disk $D$ containing all the sets, $g_n-g$ is harmonic including at
$\infty$, so $|g_n-g|\le\max_{\partial D}|g_n-g|$ there.
\end{proof}

The next proposition shows that the regularity of $\Omega$ cannot be deduced from the approximation. The Cantor-type construction in its proof is similar to the proof of Theorem~\ref{thm-sharp}: the size of subintervals is made smaller at a sequence of scales. However, this time the subintervals are made small enough so that the resulting Cantor set is polar. The approximating sequence $E_n$   is formed by pre-Cantor sets from the steps immediately preceding the small-subinterval steps.

\begin{proposition}\label{prop-irregular}
Let $\varphi\colon[1,\infty)\to (0,\infty)$ be a nondecreasing unbounded function. There exist a compact set $E\subset\R$ with $\capc E>0$ and compact $\alpha_n$-uniformly perfect sets $E_n\supset E$ with $d_n\to0$ and
$\log(1/\alpha_n) \le \varphi(t_n)$ for all $n$, such that the unbounded component
of $\widehat\C\setminus E$ is not regular.
\end{proposition}

\begin{proof} We may assume $\varphi(t)\le t$. Given numbers
$a_k\ge\log4$, let $\ell_k=\exp(-a_1-\dots-a_k)$ and
$\rho_k=\ell_k/\ell_{k-1}=e^{-a_k}$. Let $\Lambda_0=[0,1]$, obtain $\Lambda_k$
by keeping in each component of $\Lambda_{k-1}$ its two closed end intervals of
length $\ell_k$, and let $\Lambda=\bigcap_k\Lambda_k$. We need three
properties of this construction.

(a) \emph{The set $[-2,-1]\cup\Lambda_K$ is $\alpha$-uniformly perfect with
$\alpha=\frac12\min_{k\le K}\rho_k$.} Let $a\in\Lambda_K$ and $0<r\le1$, and
let $I_j$ be the component of $\Lambda_j$ containing $a$. The endpoints of
$I_j$ belong to $\Lambda_K$, and the farther of them is at distance between
$\ell_j/2$ and $\ell_j$ from $a$. If $r\le\ell_K$, the interval $I_K$ contains
a point at distance between $r/2$ and $r$ from $a$. If
$\ell_{j+1}<r\le\ell_j$ for some $j<K$, the farther endpoint of $I_{j+1}$ is at
distance between $\ell_{j+1}/2=\rho_{j+1}\ell_j/2\ge\alpha r$ and $r$ from $a$.
For $1<r\le3=\diam([-2,-1]\cup\Lambda_K)$, the farther of the points $0,1$ is at
distance between $1/2\ge\alpha r$ and $r$ from $a$. For $a\in[-2,-1]$ the
interval $[-2,-1]$ itself provides the required points, since $\alpha\le1/8$.

(b) \emph{$\Lambda\subset\Lambda_K$ and $d_H(\Lambda_K,\Lambda)\le\ell_K/2$,}
because each point of $\Lambda_K$ is within $\ell_K/2$ of an endpoint of its
component, and these endpoints belong to $\Lambda$.

(c) \emph{If $\sum_ka_k2^{-k}=\infty$, then $\Lambda$ is polar.} Let
$\nu\in\Prob(\Lambda)$. If $x$ and $y$ lie in a common component of
$\Lambda_k$, then $|x-y|\le\ell_k$. Hence
$\log(1/|x-y|)\ge\sum_ka_k\chi_k(x,y)$, where $\chi_k$ is the indicator of the
event that $x$ and $y$ lie in a common component of $\Lambda_k$. Integrating
and using the Cauchy--Schwarz inequality over the $2^k$ components $J$ of
$\Lambda_k$, we get
\[
 I(\nu)\ge 
 \iint_{}\sum_ka_k\chi_k(x,y)\,d\nu_x\,d\nu_y =  \sum_ka_k\sum_J\nu(J)^2\ge\sum_ka_k2^{-k}=\infty.
\]

Now choose integers $0=K_0<K_1<K_2<\cdots$ inductively: with
$b_n=2^{K_{n-1}}+1$, let $K_n=K_{n-1}+L_n$, where $L_n$ is so large that $\varphi(L_nb_n)\ge b_n+\log2$. Put $a_k=b_n$ for $K_{n-1}<k\le K_n$, and let
\[
 E=[-2,-1]\cup\Lambda,\qquad E_n=[-2,-1]\cup\Lambda_{K_n}.
\]
Since $(b_n)$ is increasing, (a) shows that $E_n$ is $\alpha_n$-uniformly perfect
with $\log(1/\alpha_n)=b_n+\log2$. By (b), $E\subset E_n$ and
$t_n\ge\log(2/\ell_{K_n})\ge\sum_{k\le K_n}a_k\ge L_nb_n$, so that
$\log(1/\alpha_n)\le\varphi(L_nb_n)\le\varphi(t_n)$. Each block contributes
$\sum_{K_{n-1}<k\le K_n}a_k2^{-k}\ge b_n2^{-K_{n-1}-1}\ge\frac12$, hence
$\Lambda$ is polar by (c).

A measure of finite energy does not charge a polar set. Therefore the
equilibrium measure of $E$ is that of $[-2,-1]$, and $g_E=g_{[-2,-1]}$ on
$\Omega=\widehat\C\setminus E$. The set $\Lambda$ is contained in
$\partial\Omega$, and $g_{[-2,-1]}$ is positive and continuous on $\Lambda$. Thus $\Omega$ is not
regular.
\end{proof}

\begin{proof}[Proof of Theorem~\ref{thm-green}]
If $\Omega$ is regular, then $\capc E_n\to\capc E$ by Theorem~\ref{thm-cap},
the sets $E_n$ have regular outer domains, and Lemma~\ref{lem-green} applies.
Conversely, if $g_{E_n}\to g_E$ uniformly on $\C$, then $g_E$ is continuous as
a uniform limit of the continuous functions $g_{E_n}$. Since $g_E=0$ on $\C\setminus\Omega$,
it tends to $0$ at $\partial\Omega$, and $\Omega$ is regular
\cite[Theorem 4.4.9]{Ransford}. 

Suppose that $\alpha_n\ge\alpha>0$ for all $n$. Given $a\in E$ and
$0<r<\diam E$, choose $a_n\in E_n$ with $|a-a_n|\le d_n$. For large $n$,
\eqref{eq-up} provides points of $E_n$ at distance between $\alpha r$ and $r$
from $a_n$, and they are within $d_n$ of points $w_n\in E$ with
$\alpha r-2d_n\le|w_n-a|\le r+2d_n$. A limit point of $(w_n)$ lies in $E$ at
distance between $\alpha r$ and $r$ from $a$, and the case $r=\diam E$ follows
by taking limits in $r$. Thus $E$ is $\alpha$-uniformly perfect, and $\Omega$
is regular. The last assertion of the theorem is
Proposition~\ref{prop-irregular}.
\end{proof}

\begin{remark}\label{rem-KK}
For a nonpolar limit, the proof of \cite[Theorem 3.4]{KK} contains two one-sided estimates.
The estimate $\sup_\C(g_{E_n}-g_E)\to0$ follows from the H\"older bound
\cite[Theorem 3.2]{KK}, based on Siciak~\cite{Siciak}, and the maximum
principle; it does not require regularity of $\Omega$. Letting $z\to\infty$
gives $\liminf_n\capc E_n\ge\capc E$, hence capacity convergence under the
hypothesis of \cite[Theorem 3.4]{KK}. Polar limits are already covered by
upper semicontinuity.

The reverse estimate uses~\cite[Lemma 3.5]{KK}, which cites~\cite[Lemma 2.1]{KKGreen}. The latter paper assumes all domains to be regular. The assertion in~\cite{KK} that its proof uses regularity only of
the approximating domains overlooks the relative compactness of
$\{g_E>\epsilon\}$ in $\Omega$. For every $\epsilon>0$, this property holds precisely when $\Omega$ is regular. 
The capacity conclusion of \cite[Theorem 3.4]{KK} and the NED application \cite[Theorem 4.2]{KK} remain valid.
\end{remark}

\section{Riesz capacity}\label{sec-riesz}

In this section the compact sets $E_n$ are not assumed uniformly perfect. Fix $0<s<m$. 
The Riesz capacity of a compact set $F\subset\R^m$ is
$C_s(F)=V_s(F)^{-1}$, with $1/\infty=0$. This normalization follows
Landkof~\cite{Landkof} and satisfies $C_s(rF)=r^sC_s(F)$; the length-normalized capacity used in~\cite{CL} is $C_s(F)^{1/s}$, so the
convergence statements below do not depend on the choice. The exponent $s$ is the order of the singularity of the kernel; for $m\ge3$ the Newtonian capacity corresponds to $s=m-2$, up to a constant factor. The kernel $k_s$ is strictly
positive definite, and the equilibrium measure is unique when
$C_s(F)>0$~\cite[Theorem 4.4.5]{BHS}. As in the proof of Lemma~\ref{lem-green}, if
$d_H(E_n,E)\to0$ and $C_s(E_n)\to C_s(E)>0$, then the equilibrium measures of
$E_n$ converge weakly to that of $E$.

In Theorem~\ref{thm-transfer} we have
$W_n=\sup_{a\in E_n}C_s(E_n\cap\oB(a,h_n))^{-1}$. The flexibility in the choice
of the scale $h_n\ge d_H(E_n,E)$ is convenient in examples; in particular, we choose $h_n$ to be positive even when $d_H(E_n,E)$ is zero. 

\begin{corollary}\label{cor-riesz}
Let $d_H(E_n,E)\le h_n\to0$. If there is $c>0$ such that
\begin{equation}\label{eq-cor-riesz-1}
 C_s\bigl(E_n\cap\oB(a,h_n)\bigr)\ge c\,h_n^s\qquad\text{for all $n$ and
 $a\in E_n$},    
\end{equation}
then $C_s(E_n)\to C_s(E)$.
\end{corollary}

\begin{proof}
If $C_s(E)=0$, this follows from~\eqref{eq-usc}. Otherwise let $\nu$ be the equilibrium measure of $E$. 
By the assumption~\eqref{eq-cor-riesz-1} we have $W_n\le c^{-1}h_n^{-s}$. Inequality~\eqref{eq-tail} 
gives $h_n^{-s} Q_\nu(12h_n)\to0$ as $n\to\infty$. Therefore, $W_nQ_\nu(12h_n)\le c^{-1}h_n^{-s}Q_\nu(12h_n)\to0$, and 
Corollary~\ref{cor-transfer} applies.
\end{proof}

The local capacity bound in Corollary~\ref{cor-riesz} follows from a dimensional condition. For the logarithmic kernel, uniform perfectness is such
a condition: it is equivalent to positive lower Assouad dimension~\cite[Lemma 2.1]{KLV}. For the Riesz kernel $k_s$ the dimension has to exceed $s$, and it is convenient to measure it with the Hausdorff content
\[
 \mathcal{H}^q_\infty(K)=\inf\Bigl\{\sum_i(\diam U_i)^q\colon
 K\subset\bigcup_iU_i\Bigr\}.
\]

\begin{definition}\label{def-lowerreg}
Let $0<q\le m$, $c>0$, and $h>0$. A compact set $F\subset\R^m$ is
\emph{lower $q$-regular at scales up to $h$}, with constant $c$, if
\[
 \mathcal{H}^q_\infty\bigl(F\cap\oB(a,r)\bigr)\ge c\,r^q
 \qquad\text{for all $a\in F$ and $0<r\le h$}.
\]
\end{definition}

Ahlfors $q$-regular sets are lower $q$-regular at all scales up to their
diameter. Only the scale $r=h$ is used below. The content cannot be replaced by
the Hausdorff measure $\mathcal{H}^q$: when $q<m$, a union of small balls has
infinite $\mathcal{H}^q$-measure in every ball centered on it, while its Riesz
capacity can be arbitrarily small. 

\begin{lemma}\label{lem-content}
There is $\kappa=\kappa(m)\in(0,1]$ such that for $0<s<q\le m$ every compact set
$K\subset\R^m$ satisfies
\[
 C_s(K)\ge\Bigl(1-\frac sq\Bigr)\bigl(\kappa\,\mathcal{H}^q_\infty(K)\bigr)^{s/q}.
\]
\end{lemma}

\begin{proof}
We may assume $\mathcal{H}^q_\infty(K)>0$. By Frostman's lemma  
\cite[Theorem 8.8]{Mattila},  there is a measure $\mu$ on $K$ with
$\mu(K)\ge\kappa\,\mathcal{H}^q_\infty(K)$ and $\mu(B(x,t))\le t^q$ for all
$x\in\R^m$ and $t>0$. For $\rho>0$,
\begin{align*}
 U_s^\mu(x)&=s\int_0^\infty\mu(B(x,t))\,t^{-s-1}\,dt\\
 &\le s\int_0^\rho t^{q-s-1}\,dt+s\,\mu(K)\int_\rho^\infty t^{-s-1}\,dt
 =\frac{s}{q-s}\,\rho^{q-s}+\mu(K)\rho^{-s}.
\end{align*}
The choice $\rho=\mu(K)^{1/q}$ gives
$U_s^\mu\le\frac{q}{q-s}\,\mu(K)^{1-s/q}$. Hence the probability measure
$\mu(K)^{-1}\mu$ has energy at most $\frac{q}{q-s}\,\mu(K)^{-s/q}$.
\end{proof}

\begin{theorem}\label{thm-lowerreg}
Let $d_H(E_n,E)\le h_n\to0$, and suppose that each $E_n$ is lower $q$-regular at
scales up to $h_n$ with constant $c$, where $s<q\le m$ and $c>0$ do not depend
on $n$. Then $C_s(E_n)\to C_s(E)$. More generally, if the exponents
$q_n\in(s,m]$ and the constants $c_n$ depend on $n$, the conclusion holds
provided that $C_s(E)>0$ and the equilibrium measure $\nu$ of $E$ satisfies
\begin{equation}\label{eq-lowerreg}
 h_n^{-s}Q_\nu(12h_n)=o\Bigl(\Bigl(1-\frac s{q_n}\Bigr)c_n^{s/q_n}\Bigr).
\end{equation}
\end{theorem}

\begin{proof}
Lemma~\ref{lem-content} applied to $E_n\cap\oB(a,h_n)$, together with
$\kappa^{s/q_n}\ge\kappa$, gives
\[
 W_n\le\frac{h_n^{-s}}{\kappa\,(1-s/q_n)\,c_n^{s/q_n}}.
\]
Therefore
\eqref{eq-lowerreg} implies $W_nQ_\nu(12h_n)\to0$, and
Corollary~\ref{cor-transfer} applies. When $q_n=q$ and $c_n=c$ are fixed,
\eqref{eq-lowerreg} follows from~\eqref{eq-tail} if $C_s(E)>0$, and the case
$C_s(E)=0$ follows from~\eqref{eq-usc}.
\end{proof}

Condition~\eqref{eq-lowerreg} plays the role of the condition $A_n=O(t_n)$ in
Theorem~\ref{thm-cap}: the exponent may approach $s$ and the constant may
degenerate, at a rate controlled by the concentration of $\nu$ at the scale
$h_n$.

\section*{Acknowledgments}

The first author acknowledges Saint-Petersburg State University for a research project 161726127.

Also the research of the first author was carried out within the state assignment for IAM FEB RAS 
(No.~075-00460-26-00).

\section*{AI use disclosure}

The authors used ChatGPT (OpenAI) and Claude (Anthropic) to assist with mathematical exploration and manuscript preparation. They reviewed and verified the mathematical arguments, revised the manuscript, and take full responsibility for its content.

\bibliographystyle{plain}

\bibliography{references.bib}

\end{document}